\documentclass[a4wide,11pt]{article}
\usepackage{fancybox}
\usepackage{epsfig}
\usepackage[english]{babel}
\usepackage{amssymb}
\usepackage{amsmath}
\usepackage{graphicx}
\newtheorem{counter}{Theorem}
\newtheorem{lem}{Lemma}
\newtheorem{Rmq}{Remark}
\newcommand{\A}{\mathcal{A}}
\newcommand{\C}{\mathcal{C}}
\renewcommand{\S}{\mathcal{S}}
\newcommand{\F}{\mathcal{F}}
\newcommand{\R}{\mathbb{R}}
\renewcommand{\P}{\mathbb{P}}
\newcommand{\Z}{\mathbb{Z}}
\newcommand{\N}{\mathbb{N}}
\newcommand{\M}{\mathcal{M}}
\newcommand{\U}{\mathcal{U}}
\newcommand{\V}{\mathcal{V}}
\newcommand{\Q}{\mathcal{Q}}
\renewcommand{\S}{\mathcal{S}}
\newcommand{\B}{\mathcal{B}}
\newcommand{\I}{\mathcal{I}}

\newcounter{tictac}

\def\1{\,\rlap{\mbox{\small\rm 1}}\kern.15em 1}
\def\ind#1{\1_{#1}}
\def\build#1_#2^#3{\mathrel{\mathop{\kern 0pt#1}\limits_{#2}^{#3}}}
\def\tend#1#2{\build\hbox to 12mm{\rightarrowfill}_{#1\rightarrow #2}^{ }}

\def\converge#1#2#3#4{\build\hbox to
#1mm{\rightarrowfill}_{#2\rightarrow #3}^{\hbox{\scriptsize #4}}}
\begin{document}
\thispagestyle{empty}

\vglue -2cm
\def\boxit#1#2{\vbox{\hrule\hbox{\vrule
\vbox spread#1{\vfil\hbox spread#1{\hfil#2\hfil}\vfil}%
\vrule}\hrule}}
\centerline{\large\sc CENTRE NATIONAL DE LA RECHERCHE SCIENTIFIQUE}
\smallskip
\centerline{\large\sc UNIVERSIT\'E de ROUEN }
\vspace{5cm}
\centerline{\Large\bf PUBLICATION de l'UMR 6085}
\smallskip
\centerline{\Large\bf LABORATOIRE DE MATH\'EMATIQUES RAPHA\"EL SALEM}

\vskip 7cm
\setlength{\fboxsep}{10pt}%
\centerline{
\Ovalbox{%
        \begin{Bcenter}
                \textsl{\large ON THE CENTRAL AND LOCAL LIMIT THEOREM FOR }\\
                \textsl{\large MARTINGALE DIFFERENCE SEQUENCES}\\
                \ \\
                \textbf{MOHAMED EL MACHKOURI and DALIBOR VOLN\'Y}
        \end{Bcenter}
}
}
\vfill
\centerline{\bf Document 2002-05}
\bigskip

\centerline{Universit\'e de Rouen UFR des sciences}
\centerline{Math\'ematiques, Site Colbert, UMR 6085}
\centerline{F 76821 MONT SAINT AIGNAN Cedex}
\centerline{T\'el: (33)(0) 235 14 71 00\qquad Fax: (33)(0) 232 10 37 94}
\eject

\pagestyle{myheadings}
\bibliographystyle{plain}

%%% ----------------------------------------------------------------------

\thispagestyle{empty}

\newpage

\begin{abstract}
Let $(\Omega,{\cal A},\mu)$ be a Lebesgue space and $T:\Omega\to
\Omega$  an ergodic measure preserving automorphism with positive entropy. 
We show that
there is a bounded and strictly stationary martingale difference
sequence defined on $\Omega$ with a common non-degenerate lattice
distribution satisfying the central limit theorem with an arbitrarily
slow rate of convergence and not satisfying the local limit theorem.
A similar result is established for martingale difference sequences
with densities provided the entropy is infinite. In addition, the
martingale
difference sequence may be chosen to be strongly mixing.

\vspace{1cm}\hspace{-0.7cm}
\end{abstract}

\newpage

\section{Introduction}
In this note we consider strictly stationary martingale difference
sequences $f\circ T^k$
defined on a Lebesgue probability space $(\Omega,{\cal A},\mu)$,
where $T:\Omega\to \Omega$ is an invertible measure preserving
transformation and $f:\Omega\to \mathbb R$ is measurable. 

Setting
$$ S_n(f)=\sum_{k=0}^{n-1} f\circ T^k$$
the central limit theorem we are concerned with has a formulation
that $ \mu (S_n(f)\le t\sigma \sqrt{n})$ converges to $\Phi(t)$, the
distribution function of the standard normal distribution at $t\in
\mathbb R$. This result was first obtained by de Moivre, Laplace and
Gauss, the first rigorous proof for independent, identically
distributed (i.i.d.) random variables with
finite second moment was discovered by Lyapunov around 1901. Later
Berry (1941) and Esseen (1942) showed that the rate of convergence in
this central limit theorem is of order $n^{-1/2}$ in case of existing
third moments (this rate is also best possible).
It is also well known that the central limit theorem does not imply a
local limit theorem for densities (if they exist). In 1954 Gnedenko
(see \cite{Gnedenko54}, \cite{Ibrag-Linnik} or \cite{Petrov}) solved the convergence problem of densities with the following result.
\\
\vspace{-0.2cm}
\\
{\bf Theorem A (Gnedenko, 1954)}  {\em Let $(X_n)_{n\ge 0}$ be a sequence
of i.i.d. random variables such that
$X_{0}$ has zero mean and unit variance and denote by $f_{n}$ and
$\varphi$ respectively the density function of the random variable
$n^{-1/2}(X_{0}+X_{1}+...+X_{n-1})$ and the density function of
the standard normal law. In order that
$$
\sup_{x}\vert f_{n}(x)-\varphi(x)\vert\converge{10}{n}{+\infty}{ }0
$$
it is necessary and sufficient that the density function
$f_{n_{0}}$ be bounded for some integer $n_{0}$.}\\
\vspace{-0.2cm}
\\
A similar result is valid for sums of lattice-valued random
variables. Let $(X_n)_{n\ge 0}$ be a sequence of i.i.d. random
variables having a non-degenerate distribution concentrated on $\{b+Nh:
n\in \mathbb Z\}$ where $h>0$ and $b\in \mathbb R$. $h$ is called a
maximal step if it is the largest value with this property. Suppose
the distribution has finite variance $\sigma^{2}$. Let $m=EX_0$ denote
the expectation, $S_n=X_0+...+X_{n-1}$ and 
$$ P_n(N)= \mathbb P(S_n=nb+Nh).$$
The following result has been proved by Gnedenko as well.\\
\vspace{-0.2cm}
\\
\textbf{Theorem B (Gnedenko, 1948)} {\em In order that}
$$
\sup_{N}\bigg\vert\frac{\sigma\sqrt{n}}{h}\,P_{n}(N)
-\frac{1}{\sqrt{2\pi}}\,exp\left(-\frac{1}{2}\left(\frac{nb+Nh-nm}{\sigma\sqrt{n}}\right)^{2}\right)\bigg\vert
\converge{10}{n}{+\infty}{}0
$$
{\em it is necessary and sufficient that the step $h$ should be
  maximal.}\\
\vspace{-0.2cm}
\\
Such limit theorems which deal with the local rather than with the
cumulative behaviour of the random variables are called local limit theorems (LLT).\\
Local limit theorems have been intensively studied for sums of
independent random variables and vectors together with estimates
of the rate of convergence in these theorems. In case of independent
random variables the local limit theorems can be found e.g. in \cite{Petrov},
chapter 7 (in the normal case), or in \cite{Ibrag-Linnik} (in the stable case).
In case of interval transformations the central limit theorem and
local limit theorems have been established by Rousseau-Egele \cite{Rousseau-Egele} 
and Broise \cite{Broise} using spectral theory of the Perron-Frobenius operator. The
stable local limit theorem for general transformations can be found in
Aaronson and Denker \cite{Aaronson-Denker}.

Many theorems which hold true in case of independent random variables
have been extended to martingale difference sequences, like the
central limit theorem, the law of iterated logarithm or the invariance
principle in the sense of Donsker (see \cite{Hall-Heyde}). However, there exist some
results which cannot be extended. For example, Lesigne and Voln\'y (see
\cite{Les-Vol}) have proved that the classical estimations in
large deviations inequalities for partial sums of i.i.d. random
variables cannot be attained by martingale difference sequences
even in the restricted class of ergodic and stationary sequences.
Moreover, El Machkouri and Voln\'y (see \cite{EM-Volny}) showed
that the invariance principle in the sense of Dudley may not be 
valid for martingale difference random fields.\\
In this work, we consider local and central limit theorems for
martingale difference sequences in the following setup. Let
$T:\Omega\to \Omega$ be an ergodic, measure preserving automorphism of
the Lebesgue probability space $(\Omega, {\cal A}, \mu)$ with positive
entropy.
We show (cf. Theorem $\ref{counter1}$) that one can define a stationary, bounded
martingale difference sequence with non-degenerate lattice
distribution which does not satisfy the local limit theorem as
formulated in Theorem B, but  satisfies the central limit theorem as
formulated above with an arbitrarily slow rate of convergence. In 
\cite{Peligrad-Utev}, Peligrad and Utev showed that even a very
mild mixing condition imposed on a martingale difference sequence
has a substantial impact on the limit behaviour of linear
processes generated by the sequence. In Theorem 3 we show that our examples can be
constructed with the additional property of strong (i.e. $\alpha$)
mixing.
We also give a result for martingale difference sequences with
densities (cf. Theorem $\ref{counter2}$). The rate of convergence in the central
limit theorem as in the Berry-Esseen Theorem  has been investigated for
martingale difference sequences as
well (see \cite{Hall-Heyde}).  Ibragimov (see \cite{Ibrag63}) established a Berry-Esseen type theorem for
stopped partial sums $(S_{\nu(n)}(X))_{n\geq 1}$ of bounded
martingale difference random variables $(X_{k})_{k\geq 0}$ with
the rate of convergence $n^{-1/4}$ and his estimation holds for
the whole partial sum process $(S_{n}(X))_{n\geq 1}$ if the
conditional variances of the random variables are almost surely
constant. L. Ouchti \cite{ouchti} showed that the
bounded condition in Ibragimov's result can be weakened. In 1982,
Bolthausen (see \cite{Bolth2}) has obtained the better convergence
rate of $n^{-1/2}\log\,n$ under conditions related to the behaviour
of the conditional variances.
For more recent results, see for example Haeusler \cite{Haeusler} and Jan \cite{Jan}).\\
Our results show that on the other hand, without assumptions on
the conditional variances, there is no rate of convergence for
general stationary and bounded martingale difference sequences.
\section{Main results}
Let $(X_{k})_{k\geq 0}$ be a sequence of real random variables defined
on a probability space $(\Omega, \A, \P)$ such that for any integer $k\geq 0$,
$$
E(X_{k}\vert\sigma(X_{j}\,;\,j<k))=0\quad\textrm{a.s.}
$$
Such a process is called a martingale difference sequence.
A more strict definition of the martingale difference sequence  was
introduced by Nahapetian and Petrosian
(see \cite{Nahapetian-Petrosian}): $(X_{k})_{k\geq 0}$ is a
strong martingale difference sequence if for any integer $k\geq 0$,
$$
E(X_{k}\vert\sigma(X_{j}\,;\,j\neq k))=0\quad\textrm{a.s.}
$$
Assume that the sequence $(X_{k})_{k\geq 0}$ is stationary and
ergodic, $X_{0}$ has unit variance and denote by $F$ the common
distribution function of the random variables $X_{k},\,k\geq 0$.
Recall that the Billingsley-Ibragimov central limit theorem (see
\cite{Billingsley61}, \cite{Ibrag63} or \cite{Hall-Heyde}) ensures
the convergence of the distribution functions $F_{n}$ of the
normalized partial sum process $n^{-1/2}(X_{0}+X_{1}+...+X_{n-1})$
to the standard normal distribution $\Phi$. Throughout the paper we consider the dynamical system $(\Omega, \F, \mu, T)$ where
$\Omega$ is a Lebesgue space, $\mu$ a probability measure and
$T:\Omega\to\Omega$ a bijective and bimeasurable transformation
which preserves the measure $\mu$. We denote by $\I$ the
$\sigma$-algebra of all sets $A$ in $\F$ with $TA=A$ a.s. (recall that if $\I$ is trivial then $\mu$ is said to be ergodic).
For any integer $n\geq 1$ and any zero mean random variable $f$ we denote
\begin{equation}\label{def-sum}
S_{n}(f)=\sum_{k=0}^{n-1}f\circ T^{k}
\end{equation}
and
\begin{equation}\label{def-FR}
F_{n}(f,x)=\mu(S_{n}(f)\leq x\sigma\sqrt{n}),\quad x\in\R,
\end{equation}
where $\sigma^{2}=E(f^{2})$.
\begin{counter}\label{counter1} Assume that the dynamical system
$(\Omega, \F, \mu, T)$ is ergodic and has positive entropy (cf.
\cite{Petersen} for a definition of the entropy). Let $(a_{n})_{n\geq 0}$ be a sequence of positive numbers which decreases to zero. There exists an integer
valued function $f$ in $L^{\infty}(\Omega)$ with a non-degenerate
distribution such that the following assumptions hold:
\begin{itemize}
\item the process $(f\circ T^{k})_{k\geq 0}$ is a strong martingale
  difference sequence which takes only the values -1, 0 or 1.
\item there exists an increasing sequence $(n_{k})_{k\geq 0}$ of integers such
  that for any $k\geq 0$,
\begin{equation}\label{speed-local}
\mu(S_{n_{k}}(f)=0)\geq a_{n_{k}}
\end{equation}
and
\begin{equation}\label{speed-central}
\sup_{x\in\R}\vert F_{n_{k}}(f,x)-\Phi(x)\vert\geq \frac{a_{n_{k}}}{2}.
\end{equation}
\end{itemize}
\end{counter}
\begin{Rmq} {\em By the Billingsley-Ibragimov central limit theorem for
martingale difference random variables (see \cite{Hall-Heyde}),
the sequences}
$$
\left(\mu(S_{n}(f)=0)\right)_{n\geq 1}\quad\textrm{and}\quad\left(\sup_{x\in\R}\vert F_{n}(f,x)-\Phi(x)\vert\right)_{n\geq 1}
$$
{\em converge to zero. The inequalities ($\ref{speed-local}$) and
($\ref{speed-central}$) obtained in Theorem $\ref{counter1}$ show
that this convergence can be arbitrarily slow. In particular, the
stationary process $(f\circ T^{k})_{k\geq 0}$ does not satisfy the
local limit theorem for lattice distributions (cf. Theorem B in
section 1).}
\end{Rmq}
\begin{Rmq} {\em Let $(X_{k})_{k\geq 1}$ be a sequence of bounded martingale difference random variables.
T. de la Rue \cite{de-la-Rue} proved the following
result: if there exists a positive constant $\beta$ such that for
any integer $k\geq 1$,}
\begin{equation}\label{minoration-variances}
E(X_{k+1}^{2}\vert\sigma(X_{j}\,;\,j\leq
k))\geq\beta>0\quad\textrm{a.s.}
\end{equation}
{\em then the martingale $S_{n}=X_{1}+...+X_{n}$ has a polynomial
speed of dispersion. More precisely, there exist two positive
universal constants $C$ and $\lambda$ such that for any integer
$n\geq 1$,}
$$
\sup_{t}\P(S_{n}\in I_{t})\leq Cn^{-\lambda}
$$
{\em where $I_{t}=[t-1,t+1]$ for any $t$ in $\R$. Our
counter-example (Theorem $\ref{counter1}$) shows that if the
condition ($\ref{minoration-variances}$) is not satisfied then the
speed of dispersion of the martingale $S_{n}$ can be arbitrary
slow.}
\end{Rmq}
The following result is an
analogue of Theorem $\ref{counter1}$ for sequences of martingale
difference random variables with densities.
\begin{counter}\label{counter2} 
Assume that the dynamical system $(\Omega, \F, \mu, T)$ is
ergodic and has infinite entropy. Let $(a_{n})_{n\geq 0}$ be a sequence of positive real numbers which decreases to zero. For an arbitrarily large positive
constant $L$, there exists a function $f$ in $L^{\infty}(\Omega)$
such that
\begin{itemize}
\item the random variables $\{S_{n}(f)\,;\,n\geq 1\}$ have
  densities and the density of $f$ is bounded,
\item the process $(f\circ T^{k})_{k\geq 0}$ is a strong martingale
  difference sequence,
\item there exist an increasing sequence $(n_{k})_{k\geq 0}$ of integers, a sequence $(\rho_{k})_{k\geq 0}$
  of positive real numbers which converges to zero such that for any integer $k\geq 0$,
\begin{equation}\label{oscillation-local}
\frac{1}{\rho_{k}}\mu\left(\frac{S_{n_{k}}(f)}{\sigma(f)\sqrt{n_{k}}}\in
[-\rho_{k},\rho_{k}]\right)\geq L
\end{equation}
and
\begin{equation}
\sup_{x\in\R}\vert F_{n_{k}}(f,x)-\Phi(x)\vert\geq a_{n_{k}}.
\end{equation}
\end{itemize}
\end{counter}
\begin{Rmq} {\em If the LLT holds for a stationary sequence $(g\circ T^{k})_{k\geq
0}$ of zero mean random variables then for any sequence
$(d_{k})_{k\geq 0}$ which converges to zero, we have}
$$
\frac{1}{d_{k}}\mu\left(\frac{S_{n_{k}}(g)}{\sigma(g)\sqrt{n_{k}}}\in
  [-d_{k},d_{k}]\right)\converge{10}{k}{+\infty}{ }2\varphi(0)
$$
{\em where $\varphi$ is the density of the standard normal law.
So if $L$ is sufficiently large, Inequality
($\ref{oscillation-local}$) implies that the LLT does not hold for
the strong martingale difference sequence $(f\circ T^{k})_{k\geq 0}$.}
\end{Rmq}
Now, let $(\Omega, \F, \mu)$ be a non-atomic probability space and
let $\U$ and $\V$ be two $\sigma$-algebras of $\F$. To evaluate their dependence Rosenblatt \cite{Ros} introduced the
$\alpha$-mixing coefficient defined by
$$
\alpha(\U,\V)=\sup\{\vert\mu(U\cap V)-\mu(U)\mu(V)\vert,\,U\in\U,
V\in\V\}.
$$
For any sequence $(X_{k})_{k\geq 0}$ of real random variables and
for any nonnegative integers $s$ and $t$ we denote by
$\F_{\infty}^{s}$ and $\F_{t}^{\infty}$ the $\sigma$-algebra
generated by $...,X_{s-1},X_{s}$ and $X_{t}, X_{t+1},...$
respectively. We shall use the following $\alpha$-mixing
coefficients defined for any positive integer $n$ by
\begin{equation}\label{def-mixing}
\alpha(n)=\sup_{k\geq 0}\alpha(\F_{\infty}^{k},\F_{k+n}^{\infty}).
\end{equation}
We say that the sequence $(X_{k})_{k\geq 0}$ is strongly mixing
(or $\alpha$-mixing) if $\alpha(n)$ converge to zero for $n$ going to 
infinity. For more about mixing coefficients we can refer to
Doukhan \cite{Doukhan} or Rio \cite{Rio2000}. Our last result is
the following counter-example in the topic of strongly mixing
processes.
\begin{counter}\label{counter3}
Let $(a_{n})_{n\geq 0}$ be a sequence of positive real numbers
which decreases to zero. There exists an endomorphism $T$ of
$\Omega$ and an integer-valued function $f$ in
$L^{\infty}(\Omega)$ with a non-degenerate distribution such that
\begin{itemize}
\item the process $(f\circ T^{k})_{k\geq 0}$ is a strongly mixing martingale difference sequence which takes only the values
-1, 0 and 1.
\item there exists an increasing sequence $(n_{k})_{k\geq 0}$ of integers such
that for any integer $k\geq 0$,
\begin{equation}\label{speed-local-mix}
\mu(S_{n_{k}}(f)=0)\geq a_{n_{k}}
\end{equation}
and
\begin{equation}\label{speed-central-mix}
\sup_{x\in\R}\vert F_{n_{k}}(f,x)-\Phi(x)\vert\geq
\frac{a_{n_{k}}}{2}.
\end{equation}
\end{itemize}
\end{counter}
%
% -------------------------------------------- Proofs -----------------------------------------------------
%
\section{Proofs}
\subsection{Proof of Theorem $\textbf{\ref{counter1}}$} In order to
construct the function $f$, we need the following lemma (cf. \cite{Les-Vol}).
\begin{lem}\label{lem1}
There exist two $T$-invariant sub-$\sigma$-algebras $\B$ and $\C$
of $\A$ and a $\B$-measurable function $g$ defined on $\Omega$ such that
\begin{itemize}
\item the $\sigma$-algebras $\B$ and $\C$ are independent,

\item the process $\left(g\circ T^{k}\right)_{k\geq 0}$ is a sequence
  of i.i.d. zero mean random variables which take only the
  values -1, 0 or 1,

\item the dynamical system $\left(\Omega, \C, \mu_{\C}, T\right)$
is aperiodic (that is, for any $k$ in $\Z\setminus\{0\}$ and
$\mu_{\C}$-almost all $\omega$ in $\Omega$, we have
$T^{k}\omega\neq\omega$).
\end{itemize}
Moreover, there exists $0<a\leq 1$ depending only on the
entropy of $(\Omega, \F, \mu, T)$ such that $\mu(g=\pm 1)=a/2$ and $\mu(g=0)=1-a$.
\end{lem}
The following lemma is a particular case of a result established by del Junco and Rosenblatt (\cite{delJunco-Rosenblatt}, Theorem 2.2).
\begin{lem}\label{lem2}
Consider the dynamical system $(\Sigma, \S, \nu, S)$ where $\Sigma$ is
a Lebesgue space and let us fix $\varepsilon>0$, $N$ in $\N$ and $x$ in $]0,1[$.
There exists an $\S$-measurable set $A$ such that $\nu(A)=x$ and $\nu(A\Delta
S^{-n}A)<\varepsilon\nu(A)$ for any integer $0\leq n\leq N$.
\end{lem}
For any $k\geq 0$, we fix
$\varepsilon_{k}>0,\,\,N_{k}\in\N$ and $d_{k}>0$ such that
\begin{itemize}
\item the sequence $(\varepsilon_{k})_{k\geq 0}$ decreases to zero,
\item the sequence $(N_{k})_{k\geq 0}$ increases to $+\infty$,
\item the sequence $(d_{k})_{k\geq 0}$ satisfies $\sum_{k=0}^{+\infty}
  d_{k}<1$.
\end{itemize}
Consider the $\sigma$-algebras $\B$ and $\C$ and the function $g$
given by Lemma $\ref{lem1}$. By Lemma $\ref{lem2}$, for any integer $k\geq 0$, there exists
$A_{k}$ in $\C$ such that
\begin{equation}\label{eq0}
\mu(A_{k})=d_{k}\quad\textrm{and}\quad\forall{0\leq n\leq
N_{k}}\quad\mu(A_{k}\Delta T^{-n}A_{k})<\varepsilon_{k}d_{k}.
\end{equation}
Let
\begin{displaymath}
A=\bigcup_{k=0}^{+\infty}A_{k}\in\C
\end{displaymath}
and
\begin{displaymath}
f=g\,\ind{A^{c}}\in L^{\infty}(\Omega).
\end{displaymath}
One can notice that $0<\mu(A)<1$. Moreover, the distribution of $f$ is not degenerate since
$$
\mu(f=\pm 1)=\mu(A^{c})\mu(g=\pm 1)>0.
$$
Moreover, we have for any integer $k\geq 0$,
$$
E(f\circ T^{k}\vert\F_{k})=0\quad\textrm{a.s.}
$$
where
\begin{equation}\label{tribus}
\F_{k}=\sigma(g\circ T^{j}\,;\,j\neq k)\,\vee\,\C.
\end{equation}
In particular, the stationary process $(f\circ T^{k})_{k\geq 0}$
is a strong martingale difference sequence which takes only the values -1,
0 or 1.\\
\\
{\em The Local Limit Theorem}\\
\\
Let $k\geq 0$ be a fixed integer. For any integer $1\leq n\leq N_{k}$, we have
\begin{equation}\label{inclusion1}
\bigcap_{i=0}^{n-1} T^{-i}A_{k}\subset \{S_{n}(f)=0\}
\end{equation}
Using the condition ($\ref{eq0}$), we derive
$$
\mu\left(A_{k}\backslash\bigcap_{i=0}^{n-1} T^{-i}A_{k}\right)=
\mu\left(\bigcup_{i=0}^{n-1}\left(A_{k}\backslash
T^{-i}A_{k}\right)\right)\leq N_{k}\,\varepsilon_{k}\,d_{k},
$$
hence,
$$
\mu\left(\bigcap_{i=0}^{n-1}T^{-i}A_{k}\right)\geq\mu(A_{k})-N_{k}\,\varepsilon_{k}\,d_{k}.
$$
Combining the last inequality with the assertion
($\ref{inclusion1}$), we deduce
$$
\mu(S_{n}(f)=0)\geq d_{k}(1-N_{k}\,\varepsilon_{k}).
$$
Let $(\rho_{k})_{k\geq 0}$ be a sequence of positive numbers
which converges to zero. For any $k\geq 0$ and for $\varepsilon_{k}$
sufficiently small, we can suppose that $N_{k}\,\varepsilon_{k}\leq\rho_{k}$.
Consequently, for any integer $k\geq 0$ and any $1\leq n\leq N_{k}$,
\begin{equation}\label{eq1}
\mu(S_{n}(f)=0)\geq d_{k}(1-\rho_{k}).
\end{equation}
Let $(n_{k})_{k\geq 0}$ be an increasing sequence of integers such that
$\sum_{k=0}^{+\infty}a_{n_{k}}<1/2$. For any integer $k\geq 0$ we make the
following choice:
$$
N_{k}=n_{k},\quad\rho_{k}=2^{-k-1}\quad\textrm{and}\quad
d_{k}=2a_{n_{k}}.
$$
From ($\ref{eq1}$) we deduce for any integer $k\geq 0$,
\begin{displaymath}
\mu(S_{n_{k}}(f)=0)\geq 2(1-2^{-k-1})a_{n_{k}}\geq a_{n_{k}}.
\end{displaymath}
{\em The rate of convergence in the Central Limit Theorem}\\
\\
We have
$$
a_{n_{k}}\leq\mu(S_{n_{k}}(f)=0)\leq\vert
F_{n_{k}}(f,y)-F_{n_{k}}(f,z)\vert,\quad z<0<y.
$$
From the triangular inequality it follows
$$
a_{n_{k}}\leq 2\sup_{x}\vert F_{n_{k}}(f,x)-\Phi(x)\vert+\vert\Phi(y)-\Phi(z)\vert.
$$
We have $\Phi(y)-\Phi(z)\to 0$ for $y-z$ converging to zero,
hence for any integer $k\geq 0$,
$$
\sup_{x}\vert F_{n_{k}}(f,x)-\Phi(x)\vert\geq\frac{a_{n_{k}}}{2}.
$$
The proof of Theorem $\ref{counter1}$ is complete.$\qquad\Box$
\subsection{Proof of Theorem $\textbf{\ref{counter2}}$} Since the
dynamical system $(\Omega, \F, \mu, T)$ has infinite entropy, one has the following version
of Lemma $\ref{lem1}$.
\begin{lem}\label{lem3}
There exist two $T$-invariant sub-$\sigma$-algebras $\B$ and $\C$
of $\A$ and a $\B$-measurable function $g$ defined on $\Omega$ such that
\begin{itemize}
\item the $\sigma$-algebras $\B$ and $\C$ are independent,

\item the process $\left(g\circ T^{k}\right)_{k\geq 0}$ is a sequence
  of i.i.d. random variables with common density $\ind{[-1,-1/2]}+\ind{[1/2,1]}$,

\item the dynamical system $\left(\Omega, \C, \mu_{\C}, T\right)$
is aperiodic.
\end{itemize}
\end{lem}
\textbf{Proof of Lemma $\bold{\ref{lem3}}$}. Let $(\Omega,\F,\mu)$
be a Lebesgue space and $T$ be an ergodic automorphism of
$\Omega$. W'll use the relative Sinai theorem which is
contained in the proposition $2^{'}$ in the article
\cite{Thouvenot} by Thouvenot. For any finite partition
$\M=(M_{1},...,M_{k})$ of $\Omega$, we denote by $d(\M)$ the
vector $(\mu(M_{1}),...,\mu(M_{k}))$ and by $H(M)$ the entropy of the partition.\\
\vspace{-0.3cm}
\\
\textbf{Proposition (relative Sinai theorem)} {\em Let $Q$
be a partition of $\Omega$ and let $\mathcal{P}$ be a virtual
finite partition such that $H(\mathcal{P})+H(\Q,T)\leq H(T)$. Then
there exists a partition $\mathcal{R}$ of $\Omega$ which satisfies
\begin{itemize}
\item the sequence $\{T^{i}\mathcal{R}\}_{i\in\Z}$ is independent,
\item $d(\mathcal{R})=d(\mathcal{P})$,
\item $\bigvee_{-\infty}^{+\infty}T^{i}\mathcal{R}$ and $\bigvee_{-\infty}^{+\infty}T^{i}\Q$ are
independent.
\end{itemize}}
\hspace{-0.77cm} Recall that the entropy of the ergodic dynamical
system $(\Omega,\F,\mu,T)$ is assumed to be infinite. Using the
relative Sinai theorem we construct by induction a sequence
$\{\Q_{s}\}_{s\geq 0}$ of finite partitions such that for any
integer $s\geq 0$,
\begin{itemize}
\item the sequence $\{T^{i}\Q_{s}\}_{i\in\Z}$ is independent,
\item $\Q_{s}=\{A_{s},A_{s}^{c}\}$ with $\mu(A_{s})=1/2$, hence,
$H(\Q_{s})=1$
\end{itemize}
and such that the $\sigma$-algebras
$\{\bigvee_{-\infty}^{+\infty}T^{i}\Q_{s}\}_{s\geq 0}$ are
mutually independent. Denote by $\B_{0}$ the $\sigma$-algebra
generated by the partitions $\{\Q_{s},\,s\geq 1\}$ and remark that
the sequence $\{T^{i}\B_{0}\}_{i\in\Z}$ is independent. Now 
consider the following two invariant sub-$\sigma$-algebras of $\F$
$$
\C=\bigvee_{-\infty}^{+\infty}T^{i}\mathcal{Q}_{0}
\quad\textrm{and}\quad \B=\bigvee_{-\infty}^{+\infty}T^{i}\B_{0}.
$$
By construction, the $\sigma$-algebras $\B$ and $\C$ are
independent. Moreover, there exists a $\B_{0}$-measurable function
$\psi_{1}$ defined on $\Omega$ to $[0,1]$ such that for any
integer $s\geq 1$,
$$
A_{s}=\psi_{1}^{-1}\left(\{x\in[0,1]\,;\,2^{s-1}x-[2^{s-1}x]<1/2\}\right)
$$
where $[\,.\,]$ denotes the integer part function. Let $\psi_{2}$
be the function defined on $[0,1]$ to $[-1,-1/2[\cup[1/2,1]$ by
$$
\psi_{2}(x)=\left\{\begin{array}{ll}
             x-1 &\textrm{if $\quad 0\leq x<1/2$}\\
             x &\textrm{if $\quad 1/2\leq x\leq 1$}.
             \end{array}
             \right.
$$
Finally the function $g=\psi_{2}\circ\psi_{1}$ satisfies the
required properties.
The proof of Lemma $\ref{lem3}$ is complete.$\qquad\Box$\\
\vspace{-0.2cm}
\\
Let $(p_{k})_{k\geq 0}$ be a sequence of positive real numbers such that
\begin{equation}\label{speed-p-k}
\sum_{k=0}^{+\infty}p_{k}=1
\end{equation}
and let $(N_{k})_{k\geq 0}$ be an increasing sequence of positive integers
such that the greatest common divisor of all $N_{k}$ be 1.\\
By the multiple Rokhlin tower theorem of Alpern (see
\cite{Alpern}), there exist $\C$-measurables sets $(F_{k})_{k\geq
0}$ with $\mu(F_{k})=p_{k}/N_{k}$ for any integer $k\geq 0$ such
that $\{T^{i}F_{k}, 0\leq i\leq N_{k}-1, k\in\N\}$ is a partition
of $\Omega$. Let $(d_{k})_{k\geq 0}$ be a sequence of positive
real numbers which converges to zero and let us denote
\begin{equation}\label{notation1}
G_{k}=\bigcup_{i=0}^{N_{k}-1}T^{i}F_{k},
\end{equation}
\begin{equation}\label{notation2}
h=\sum_{k=0}^{+\infty}d_{k}\ind{G_{k}},
\end{equation}
and
\begin{equation}\label{notation3}
f=g\,h\in L^{\infty}(\Omega)
\end{equation}
where $g$ is the function given by Lemma $\ref{lem3}$. Considering
the $\sigma$-algebras defined by ($\ref{tribus}$) we derive (as
in the proof of Theorem $\ref{counter1}$) that the stationary
process $(f\circ T^{k})_{k\geq 0}$ is a strong martingale
difference sequence.\\
\\
{\em The densities of the partial sums $\{S_{n}(f),\,n\geq 1\}$}\\
\\
Let us show that the random variables $\{S_{n}(f)\,;\,n\geq 1\}$
defined by ($\ref{def-sum}$) have densities. Let $\mathcal{P}$ be
the partition $\{G_{0},G_{1},...\}$ of $\Omega$ and let $n$ be a
fixed positive integer. Consider the following partition of
$\Omega$
$$
\mathcal{P}_{n}=\bigvee_{j=0}^{n-1}T^{-j}\mathcal{P},
$$
let $A$ in $\mathcal{P}_{n}$ be fixed and let
$r_{0}(A),r_{1}(A),...,r_{n-1}(A)$ be nonnegative integers such
that
$$
A=\bigcap_{j=0}^{n-1}T^{-j}G_{r_{j}(A)}.
$$
For any real $x$, define
$$
F_{A}(S_{n},x)=\mu(A\cap\{S_{n}(f)\leq x\})\quad\textrm{and}\quad
F_{n}(x)=\mu(S_{n}(f)\leq x).
$$
Let $0\leq j<n$ be fixed and let $d_{j}(A)$ denote $d_{r_{j}(A)}$.
For any $\omega$ in $A$ we can check that
$f(T^{j}\omega)=d_{j}(A)\,g(T^{j}\omega)$, hence for any real
$x$,
$$
A\cap\{S_{n}(f)\leq
x\}=A\cap\bigg\{\sum_{j=0}^{n-1}d_{j}(A)\,g\circ T^{j}\leq x\bigg\}.
$$
By the independence of the $\sigma$-algebras $\B$ and $\C$ we
thus get that for any real number $x$,
$$
F_{A}(S_{n},x)=\mu(A)\,\mu\left(\sum_{j=0}^{n-1}d_{j}(A)\,g\circ T^{j}\leq x\right).
$$
The distribution function of $g$ is differentiable at all
$x\notin\{-1,-1/2,1/2,1\}$, hence, the function $F_{A}(S_{n},.)$
is differentiable at all real $x$ except a finite set.\\
Let $x<0$ be a fixed real number. Since the sequence
$(d_{k})_{k\geq 0}$ converges to zero, there is only a finite
number of sets $A$ in $\mathcal{P}_{n}$ such that $\mu(A)>0$ and
$\sum_{j=0}^{n-1}d_{j}(A)>-x$. Consequently, the sum
$$
F_{n}(x)=\sum_{A\in\mathcal{P}_{n}}F_{A}(S_{n},x)
$$
contains only finitely many nonzero terms. So the distribution
function $F_{n}$ of $S_{n}(f)$ is differentiable on
$\R_{-}^{\ast}$ except a countable set. By the symmetry of the
density of $g$ and the independence of the process $(g\circ
T^{k})_{k\geq 0}$, we deduce the differentiability of $F_{n}$ at
all points $x$ in $\R$ except a countable set. Hence for any
positive integer $n$, the random variable $S_{n}(f)$ has a
density. Now we will show that the density of the random
variable $f$ (or $S_{1}(f)$) is bounded. Consider the distribution
function $F_{1}$ of the random variable $f$. Using the
independence of the $\sigma$-algebras $\B$ and $\C$ we have for
any real $x$
\begin{align*}
F_{1}(x)&=\sum_{k=0}^{+\infty}\mu(G_{k}\cap\{d_{k}g\leq x\})\\
&=\sum_{k=0}^{+\infty}\mu(G_{k})\,\mu(d_{k}g\leq x)\\
&=\sum_{k=0}^{+\infty}p_{k}\,\mu(d_{k}g\leq x).
\end{align*}
Moreover, for any integer $k\geq 0$ and any real $x$,
$$
\mu(d_{k}g\leq x)=\left\{\begin{array}{lll}
                        1 &\textrm{if} & x\geq d_{k}\\
                        \frac{x}{d_{k}}&\textrm{if}&\frac{d_{k}}{2}\leq x<d_{k}\\
                        \frac{1}{2} &\textrm{if} &-\frac{d_{k}}{2}\leq x<\frac{d_{k}}{2}\\
                        \frac{x+d_{k}}{d_{k}}&\textrm{if}&-d_{k}\leq x<-\frac{d_{k}}{2}\\
                        0 &\textrm{if} & x\leq -d_{k}.
                        \end{array}
                        \right.
$$
Let $L_{1}$ and $L_{2}$ be two positive constants such that
$L_{2}\gg L_{1}$. In the sequel, for any integer $k\geq 0$ we put
$$
p_{k}=\frac{\sqrt{2}-1}{\sqrt{2}}2^{-k/2}
$$
and
$$
d_{k}=\left\{\begin{array}{ll}
             p_{k}/L_{1} &\textrm{if $k$ is even}\\
             p_{k}/L_{2} &\textrm{if $k$ is odd}.
             \end{array}
             \right.
$$
Thus, the condition ($\ref{speed-p-k}$) still holds and the
function $f$ defined by ($\ref{notation3}$) depends on the
constants $L_{1}$ and $L_{2}$. If we put
$$
c_{1}=\sum_{j\geq 0}p_{2j}^{3}\quad\textrm{and}\quad
c_{2}=\sum_{j\geq 0}p_{2j+1}^{3}
$$
then the variance of the function $f$ is given by
\begin{equation}\label{def-variance-f}
\sigma^{2}(f)=\frac{7}{12}\left(\frac{c_{1}}{L_{1}^{2}}+\frac{c_{2}}{L_{2}^{2}}\right).
\end{equation}
If $x\geq d_{0}$ then $x\geq d_{k}$ for any integer $k\geq 0$,
hence
$$
F_{1}(x)=\sum_{j=0}^{+\infty}p_{j}=1.
$$
If $0<x<d_{0}$ then there exists a unique odd integer $k=k(x)$
such that $d_{k+2}=d_{k}/2\leq x<d_{k}$ and there exists a unique
even integer $l=l(x)$ such that $d_{l+2}=d_{l}/2\leq x<d_{l}$. So
we have
$$
\sum_{j=0}^{+\infty}p_{2j+1}\,\mu(d_{2j+1}g\leq
x)=\sum_{\substack{j\leq k-1 \\ \textrm{$j$ odd}}
}\frac{p_{j}}{2}+x\frac{p_{k}}{d_{k}}+\sum_{\substack{j\geq k+1 \\
\textrm{$j$ odd}}}p_{j}
$$
and
$$
\sum_{j=0}^{+\infty}p_{2j}\,\mu(d_{2j}g\leq
x)=\sum_{\substack{j\leq l-1 \\ \textrm{$j$ even}}
}\frac{p_{j}}{2}+x\frac{p_{l}}{d_{l}}+\sum_{\substack{j\geq l+1 \\
\textrm{$j$ even}}}p_{j}.
$$
Consequently the distribution function $F_{1}$ of $f$ is
differentiable at all points $x>0$ which do not belong to the set
$\{d_{k}\,;\,k\geq 0\}$. Since $F_{1}$ is symmetric we obtain its
differentiability at all points of the real line except a
countable set and the density $F_{1}^{'}$ of $f$ is bounded by
$L_{1}+L_{2}$.\\
\\
{\em The Local Limit Theorem}\\
\\
For any integer $n\geq 1$ we denote
$$
b_{n}=\mu\left(n^{-1/2}S_{n}(g)\in [-1,1]\right)
$$
From the central limit theorem it follows that there exists a
positive real number $b$ such that $b_{n}\geq b>0$ for all
sufficiently large $n$. Moreover, there exists an increasing
sequence $(n_{k})_{k\geq 0}$ of integers such that for any $k\geq
0$
\begin{equation}\label{condition-rho-k-a-k}
\rho_{k}\geq a_{n_{k}}
\end{equation}
where $\rho_{k}=d_{k}/\sigma(f)$. Let $k\geq 0$ be a fixed
integer, choose $N_{k}$ such that $N_{k}\geq 2n_{k}$ and denote
$$
\widetilde{G}_{k}=\bigcup_{i=0}^{N_{k}-n_{k}}T^{i}F_{k}\subset G_{k}
$$
and
$$
E_{k}=\{\omega\in\widetilde{G}_{k}\,;\,n_{k}^{-1/2}\sum_{i=0}^{n_{k}-1}f(T^{i}\omega)\in
[-d_{k},d_{k}]\}.
$$
Let $i$ in $\{0,...,n_{k}-1\}$ be fixed. Since
$T^{i}\widetilde{G}_{k}\subset G_{k}$, we have
$h(T^{i}\omega)=d_{k}$ for any $\omega$ in $\widetilde{G}_{k}$;
using ($\ref{notation3}$) we deduce
\begin{align*}
E_{k}&=\{\omega\in\widetilde{G}_{k}\,;\,n_{k}^{-1/2}\sum_{i=0}^{n_{k}-1}g(T^{i}\omega)
\,h(T^{i}\omega)\in [-d_{k},d_{k}]\}\\
&=\{\omega\in\widetilde{G}_{k}\,;\,n_{k}^{-1/2}\sum_{i=0}^{n_{k}-1}d_{k}
\,g(T^{i}\omega)\in [-d_{k},d_{k}]\}\\
&=\{\omega\in\widetilde{G}_{k}\,;\,n_{k}^{-1/2}\sum_{i=0}^{n_{k}-1}g(T^{i}\omega)\in
[-1,1]\}.
\end{align*}
By the independence of $\B$ and $\C$ it follows that
$$
\mu\left(n_{k}^{-1/2}S_{n_{k}}(f)\in [-d_{k},d_{k}]\right)
\geq\mu(E_{k})=\mu(\widetilde{G}_{k})\,b_{n_{k}}\geq\mu(\widetilde{G}_{k})\,b.
$$
Moreover
$$
\mu(\widetilde{G}_{k})=p_{k}(1-\frac{n_{k}}{N_{k}})\geq\frac{p_{k}}{2}.
$$
Hence, we derive
$$
\frac{1}{d_{k}}\mu\left(n_{k}^{-1/2}S_{n_{k}}(f)\in
[-d_{k},d_{k}]\right)\geq\frac{bp_{k}}{2d_{k}}.
$$
From th equality $\rho_{k}=d_{k}/\sigma(f)$ and from 
($\ref{def-variance-f}$) it follows that
$$
\frac{1}{\rho_{k}}\mu\left(\frac{S_{n_{k}}(f)}{\sigma(f)\sqrt{n_{k}}}\in
[-\rho_{k},\rho_{k}]\right)\geq\frac{bp_{k}}{2d_{k}}\sigma(f)
=\frac{\sqrt{7}bp_{k}}{4\sqrt{3}d_{k}}\left(\frac{c_{1}}{L_{1}^{2}}+\frac{c_{2}}{L_{2}^{2}}\right)^{1/2}.
$$
Consequently, if $k$ is odd then
$$
\frac{1}{\rho_{k}}\mu\left(\frac{S_{n_{k}}(f)}{\sigma(f)\sqrt{n_{k}}}\in
[-\rho_{k},\rho_{k}]\right)\geq\frac{bp_{k}}{2d_{k}}\sigma(f)
=\frac{\sqrt{7}b}{4\sqrt{3}}\left(\frac{c_{1}L_{2}^{2}}{L_{1}^{2}}+c_{2}\right)^{1/2}.
$$
Finally choosing $L_{2}$ sufficiently large we derive that for
$k$ odd
\begin{equation}\label{speed-local-density}
\frac{1}{\rho_{k}}\mu\left(\frac{S_{n_{k}}(f)}{\sigma(f)\sqrt{n_{k}}}\in
[-\rho_{k},\rho_{k}]\right)\geq L.
\end{equation}
As a consequence the strong martingale difference sequence
$(f\circ T^{k})_{k\geq 0}$ does not satisfy the local limit
theorem for densities.\\
\\
{\em The rate of convergence in the Central Limit Theorem}\\
\\
Now we are going to prove the last part of Theorem
$\ref{counter2}$. Recall that $\Phi$ and $\varphi$ are
respectively the distribution function and the density function of
the standard normal law. Let $k\geq 0$ be a fixed odd integer.
Using ($\ref{speed-local-density}$) we obtain that
\begin{align*}
L\rho_{k}&\leq\mu\left(\frac{S_{n_{k}}(f)}{\sigma(f)\sqrt{n_{k}}}\in
  [-\rho_{k},\rho_{k}]\right)\\
&=\vert F_{n_{k}}(f,\rho_{k})-F_{n_{k}}(f,-\rho_{k})\vert\\
&\leq 2\sup_{x}\vert F_{n_{k}}(f,x)-\Phi(x)\vert+\vert\Phi(\rho_{k})-\Phi(-\rho_{k})\vert\\
&\leq 2\sup_{x}\vert
F_{n_{k}}(f,x)-\Phi(x)\vert+2\rho_{k}\varphi(0),
\end{align*}
hence
$$
\sup_{x}\vert F_{n_{k}}(f,x)-\Phi(x)\vert\geq
(L/2-\varphi(0))\rho_{k}.
$$
Putting $L$ sufficiently large and using inequality
($\ref{condition-rho-k-a-k}$) we derive that for any odd integer
$k\geq 0$,
$$
\sup_{x}\vert F_{n_{k}}(f,x)-\Phi(x)\vert\geq a_{n_{k}}.
$$
The proof of Theorem $\ref{counter2}$ is complete.$\qquad\Box$
\subsection{Proof of Theorem $\textbf{\ref{counter3}}$}
Consider the non-atomic probability space $(\Omega, \A, \mu)$. By induction we shall construct a sequence $(p_{k})_{k\geq 1}$ of
positive real numbers with $\sum_{k\geq 1}p_{k}=1$, positive
integers $(N_{k})_{k\geq 1}$ and $(n_{k})_{k\geq 1}$, measurable
sets $(F_{k})_{k\geq 1}$ and a bijective bimeasurable
transformation $T:\Omega\to\Omega$ which preserves the measure
$\mu$ such that the sets
$T^{j}F_{k},\,j=0,...,N_{k}-1,\,k=1,2,...$ form a partition of
$\Omega$ in the way that
$$
\mu\left(\bigcup_{j=0}^{N_{k}-1}T^{j}F_{k}\right)=p_{k},\,\,k\geq
1
$$
and such that for any $k\geq 1$
\begin{equation}\label{star}
p_{k}\geq 4a_{n_{k}}\quad\textrm{and}\quad N_{k}\geq 4 n_{k}^{2}.
\end{equation}
Assume that this construction is achieved. Let $k\geq 1$ be a
fixed integer and let $\overline{F}_{k}\subset F_{k}$ be a
measurable set such that $\mu(\overline{F}_{k})=\mu(F_{k})/2$ and
define
$$
A_{k}=\bigcup_{j=0}^{N_{k}-n_{k}}T^{j}\overline{F}_{k}.
$$
Let $(g\circ T^{k})_{k\geq 0}$ be a sequence of independent random
variables independent of the $\sigma$-algebra
$\bigvee_{i=-\infty}^{+\infty}T^{i}\sigma(\overline{F}_{k},F_{k}\backslash\overline{F}_{k},
k\geq 1)$ and such that $\mu(g=\pm 1)=1/2$. Denote by
\begin{equation}\label{def-f}
f=g\ind{A^{c}}\quad\textrm{where}\quad A=\bigcup_{k\geq 1}A_{k}.
\end{equation}
As done in the proof of Theorem $\ref{counter2}$ we can check
that the process $(f\circ T^{k})_{k\geq 0}$ is a strong martingale
difference sequence. Let $k\geq 1$ be a fixed integer. By
construction we have
\begin{align*}
\mu(A_{k}\backslash\cap_{i=0}^{n_{k}-1}T^{-i}A_{k})&\leq\sum_{i=0}^{n_{k}-1}\mu(A_{k}\backslash
T^{-i}A_{k})\\
&=\sum_{i=0}^{n_{k}-1}i\times p_{k}/2N_{k}\\
&\leq\frac{n_{k}^{2}p_{k}}{2N_{k}}
\end{align*}
So using ($\ref{star}$) we derive
\begin{align*}
\mu(\cap_{i=0}^{n_{k}-1}T^{-i}A_{k})
&\geq\mu(A_{k})-\frac{n_{k}^{2}p_{k}}{2N_{k}}\\
&=\frac{(N_{k}-n_{k}+1)p_{k}}{2N_{k}}-\frac{n_{k}^{2}p_{k}}{2N_{k}}\\
&\geq\frac{p_{k}}{2}-\frac{n_{k}^{2}p_{k}}{N_{k}}\\
&=p_{k}\left(\frac{1}{2}-\frac{n_{k}^{2}}{N_{k}}\right)\\
&\geq\frac{p_{k}}{4}\geq a_{n_{k}}.
\end{align*}
Since $\{S_{n_{k}}(f)=0\}\supset\cap_{i=0}^{n_{k}-1}T^{-i}A_{k}$
we obtain inequality ($\ref{speed-local-mix}$). The proof of
Inequality ($\ref{speed-central-mix}$) is obtained from Inequality
($\ref{speed-local-mix}$) as done in section 3.1.2, so it is left
to the reader. Now, using induction, we will construct the
sequences $(p_{k})_{k}$, $(N_{k})_{k}$, $(n_{k})_{k}$,
$(F_{k})_{k}$ and the transformation $T$ such that the condition
($\ref{star}$) holds and we have to show that the martingale
difference sequence
$(f\circ T^{k})_{k\geq 0}$ is also strongly mixing.\\
%
%---------------------------------------------- recurrence ---------------------------------------------------------------
%
Let $k\geq 2$ be a fixed integer, let $p_{1},...,p_{k-1}>0$ be
given with $\sum_{i=1}^{k-1}p_{i}<1$ and denote
$p_{k}^{'}=1-\sum_{i=1}^{k-1}p_{i}$. Let $n_{1}<n_{2}<...<n_{k-1}$
and $N_{1}<N_{2}<...<N_{k-1}$ be sufficiently large integers such
that the condition ($\ref{star}$) holds. There exist $k$ sets
$F_{l},\,l=1,...,k-1,\,F_{k}^{'}$ and a bijective bimeasurable
measure-preserving transformation $T_{k}:\Omega\to\Omega$ such
that the sets $T_{k}^{j}F_{l},\,j=0,...,N_{l}-1,\,l=1,...,k-1,\,
T_{k}^{j}F_{k}^{'},\,j=0,...,N_{k}-1$ form a partition of $\Omega$
and
$$
\mu\left(\bigcup_{j=0}^{N_{k}-1}T_{k}^{j}F_{k}^{'}\right)=p_{k}^{'},\quad
\mu\left(\bigcup_{j=0}^{N_{l}-1}T_{k}^{j}F_{l}\right)=p_{l},\,l=1,...,k-1.
$$
On the sets
$T_{k}^{N_{l}-1}F_{l},\,l=1,...,k-1,\,T_{k}^{N_{k}-1}F_{k}^{'}$
the transformation $T_{k}$ acts in the way that
$$
T_{k}\left(T_{k}^{N_{k}-1}F_{k}^{'}\cup\bigcup_{l=1}^{k-1}T_{k}^{N_{l}-1}F_{l}\right)=F_{k}^{'}\cup\bigcup_{l=1}^{k-1}F_{l}
$$
and for any $F$ in
$\{T_{k}^{N_{k}-1}F_{k}^{'},T_{k}^{N_{1}-1}F_{1},...,T_{k}^{N_{k-1}-1}F_{k-1}\}$
$F$ is mapped into $F_{k}^{'},\,F_{1},...,F_{k-1}$ with
probabilities
$$
\mu(F_{k}^{'}\vert TF)=
\frac{p_{k}^{'}/N_{k}}{p_{k}^{'}/N_{k}+\sum_{j=1}^{k-1}p_{j}/N_{j}}
$$
and
$$
\mu(F_{l}\vert TF)=
\frac{p_{l}/N_{l}}{p_{k}^{'}/N_{k}+\sum_{j=1}^{k-1}p_{j}/N_{j}},\,\,l=1,...,k-1.
$$
Let us construct the transformation $T_{k+1}$. We choose an
arbitrarily small $\delta_{k}>0$ and define
$p_{k}=(1-\delta_{k})p_{k}^{'}$ and
$p_{k+1}^{'}=\delta_{k}p_{k}^{'}$. We choose also a positive
integer $N_{k+1}$ sufficiently large such that the condition
($\ref{star}$) holds. Then there exist sets $F_{k}\subset
F_{k}^{'}$ and $F_{k+1}^{'}\subset F_{k}^{'}$ and an automorphism
$T_{k+1}:\Omega\to\Omega$ such that
\begin{itemize}
\item $T_{k+1}=T_{k}$ on the set
$\bigcup_{l=1}^{k}\bigcup_{j=0}^{N_{l}-2}T_{k}^{j}F_{l}$\\
\item the sets $T_{k+1}^{j}F_{l},\,j=0,...,N_{l}-1,\,l=1,...,k$
(which are equal to the sets
$T_{k}^{j}F_{l},\,j=0,..,N_{l}-1,\,l=1,...,k$) and
$T_{k+1}^{j}F_{k+1}^{'},\,j=0,...,N_{k+1}-1$ form a partition of $\Omega$.\\
\item $\mu\left(\cup_{j=0}^{N_{k}-1}T_{k+1}^{j}F_{k}\right)=p_{k}$ and
$\mu\left(\cup_{j=0}^{N_{k+1}-1}T_{k+1}^{j}F_{k+1}^{'}\right)=p_{k+1}^{'}$.
\end{itemize}
On the sets
$T_{k+1}^{N_{k+1}-1}F_{k+1}^{'},\,T_{k+1}^{N_{l}-1}F_{l},\,l=1,...,k$,
we define $T_{k+1}$ so that
$$
T_{k+1}\left(T_{k+1}^{N_{k+1}-1}F_{k+1}^{'}\cup\bigcup_{l=1}^{k}T_{k+1}^{N_{l}-1}F_{l}\right)=F_{k+1}^{'}\cup\bigcup_{l=1}^{k}F_{l}
$$
and for any $F$ in
$\{T_{k+1}^{N_{k+1}-1}F_{k+1}^{'},T_{k+1}^{N_{1}-1}F_{1},...,T_{k+1}^{N_{k}-1}F_{k}\}$,
$F$ is mapped into $F_{k+1}^{'},\,F_{1},...,F_{k}$ with
probabilities
$$
\mu(F_{k+1}^{'}\vert TF)=
\frac{p_{k+1}^{'}/N_{k+1}}{p_{k+1}^{'}/N_{k+1}+\sum_{j=1}^{k}p_{j}/N_{j}}
$$
and
$$
\mu(F_{l}\vert
TF)=\frac{p_{l}/N_{l}}{p_{k+1}^{'}/N_{k+1}+\sum_{j=1}^{k}p_{j}/N_{j}},\,\,l=1,...,k.
$$
Consider the sets
$G_{l}=\cup_{j=0}^{N_{l}-1}T_{k}^{j}F_{l},\,l=1,...,k-1$,
$G_{k}^{'}=\cup_{j=0}^{N_{k}-1}T_{k}^{j}F_{k}^{'}$ and
$G_{k+1}^{'}=\cup_{j=0}^{N_{k+1}-1}T_{k+1}^{j}F_{k+1}^{'}$ and
denote by $\Delta_{k}$ the set $\cup_{l=1}^{k}G_{l}$. In order to clarify the construction we present the following picture which shows the $k$-st and $(k+1)$-st partitions of $\Omega$.
%%%%%%%%%%%%%%%%%%%%%%%%%%%% figure %%%%%%%%%%%%%%%%%%%%%%%%%%%%
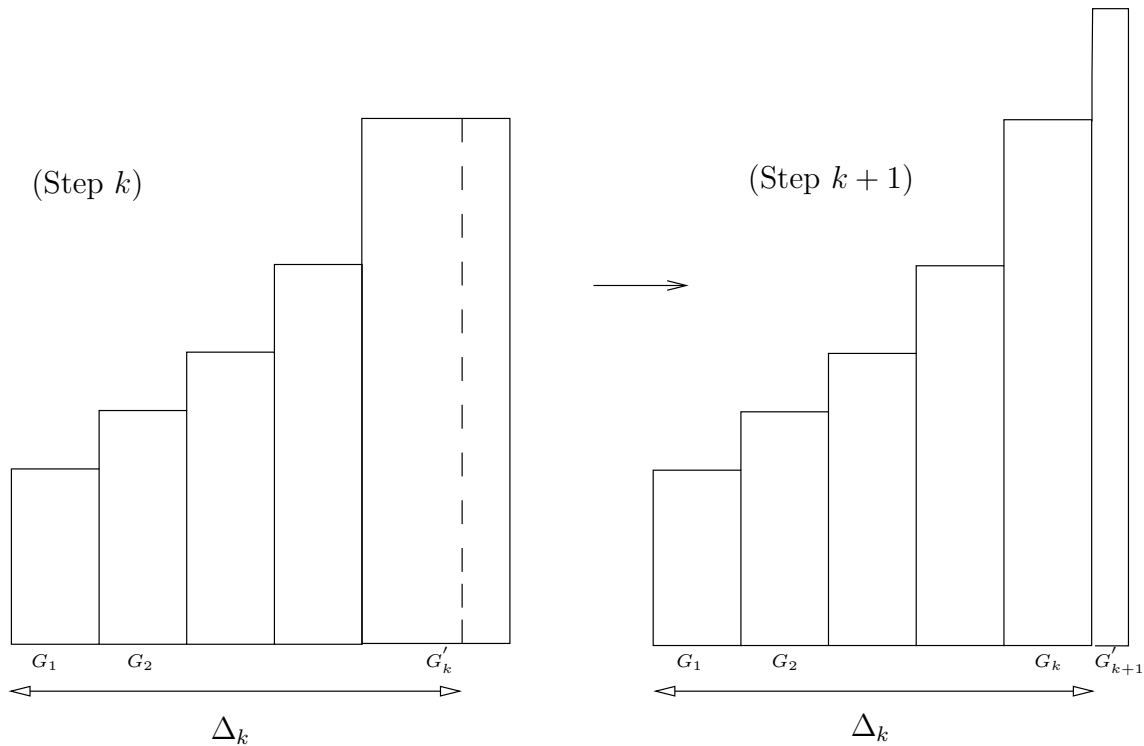
\begin{center}
\begin{picture}(0,0)%
\includegraphics{tourtronquee.pstex}%
\end{picture}%
\setlength{\unitlength}{3947sp}%
\begingroup\makeatletter\ifx\SetFigFont\undefined%
\gdef\SetFigFont#1#2#3#4#5{%
  \reset@font\fontsize{#1}{#2pt}%
  \fontfamily{#3}\fontseries{#4}\fontshape{#5}%
  \selectfont}%
\fi\endgroup%
\begin{picture}(7059,4602)(1189,-3811)
\put(1351,-376){\makebox(0,0)[lb]{\smash{\SetFigFont{12}{14.4}{\rmdefault}{\mddefault}{\updefault}(Step $k$)}}}
\put(5851,-346){\makebox(0,0)[lb]{\smash{\SetFigFont{12}{14.4}{\rmdefault}{\mddefault}{\updefault}(Step $k+1$)}}}
\put(6496,-3796){\makebox(0,0)[lb]{\smash{\SetFigFont{12}{14.4}{\rmdefault}{\mddefault}{\updefault}$\Delta_k$}}}
\put(2476,-3811){\makebox(0,0)[lb]{\smash{\SetFigFont{12}{14.4}{\rmdefault}{\mddefault}{\updefault}$\Delta_k$}}}
\put(1351,-3361){\makebox(0,0)[lb]{\smash{\SetFigFont{7}{8.4}{\rmdefault}{\mddefault}{\updefault}$G_{1}$}}}
\put(1951,-3361){\makebox(0,0)[lb]{\smash{\SetFigFont{7}{8.4}{\rmdefault}{\mddefault}{\updefault}$G_{2}$}}}
\put(3826,-3361){\makebox(0,0)[lb]{\smash{\SetFigFont{7}{8.4}{\rmdefault}{\mddefault}{\updefault}$G_{k}^{'}$}}}
\put(5401,-3361){\makebox(0,0)[lb]{\smash{\SetFigFont{7}{8.4}{\rmdefault}{\mddefault}{\updefault}$G_{1}$}}}
\put(6001,-3361){\makebox(0,0)[lb]{\smash{\SetFigFont{7}{8.4}{\rmdefault}{\mddefault}{\updefault}$G_{2}$}}}
\put(7651,-3361){\makebox(0,0)[lb]{\smash{\SetFigFont{7}{8.4}{\rmdefault}{\mddefault}{\updefault}$G_{k}$}}}
\put(8026,-3361){\makebox(0,0)[lb]{\smash{\SetFigFont{7}{8.4}{\rmdefault}{\mddefault}{\updefault}$G_{k+1}^{'}$}}}
\end{picture}

\end{center}
%%%%%%%%%%%%%%%%%%%%%%%%%%%%%%%%%%%%%%%%%%%%%%%%%%%%%%%%%%%%%%%%
Moreover, by choosing $\delta_{k}$ sufficiently small we are able
to make the measure of the set
$\{\omega\in\Omega\,,\,T_{k}\omega\neq T_{k+1}\omega\}$ to be as
small as we wish. For almost every $\omega$ in $\Omega$ there
then exists $T\omega$ in $\Omega$ such that $T_{k}\omega=T\omega$
for all $k$ sufficiently large. The transformation
$T:\Omega\to\Omega$ acts in the way that
$$
T\left(\bigcup_{l\geq 1}T^{N_{l}-1}F_{l}\right)=\bigcup_{l\geq 1
}F_{l}
$$
and for any $F$ in $\{T^{N_{1}-1}F_{1},T^{N_{2}-1}F_{2},...\}$
$F$ is mapped into $F_{1},F_{2},...$ with probabilities
$$
\frac{p_{l}/N_{l}}{\sum_{j\geq 1}p_{j}/N_{j}},\,\,l\geq 1.
$$
By $\S$ we denote the family of the sets
$\{T^{j}F_{l},\,j=0,...,N_{l}-1,\,l=1,2,...\}$ and by $\S_{k}$ we
denote the family
$\{T_{k}^{j}F_{l},\,j=0,...,N_{l}-1,\,l=1,...,k-1\}\cup\{T_{k}^{j}F_{k}^{'},\,j=0,...,N_{k}-1\}$.
The transformation $T$ then defines a Markov chain $(\xi_{i})_{i}$
with the state space $\S$ (in the way that if $T^{i}\omega\in
T^{j}F_{l}$ then $\xi_{i}(\omega)=T^{j}F_{l}$) and the
transformation $T_{k}$ defines a Markov chain
$(\xi_{i}^{(k)})_{i}$ with the state space $\S_{k}$ (analogically). 
We choose the numbers $(N_{l})_{l\geq 1}$ so that for any
$k\geq 2$, the greatest common divisor of $\{N_{1},...,N_{k}\}$ is
1. Because $T$ and $T_{k}$ are automorphisms of $\Omega$, the
chains $(\xi_{i})_{i}$ and $(\xi_{i}^{(k)})_{i}$ are stationary
and due to the choice of the numbers $(N_{l})_{l\geq 1}$ they
are aperiodic and irreducible. For any $k$ the state space
$\S_{k}$ is finite, hence for any $\varepsilon>0$ there exists
$m\in\N$ (cf. Billingsley \cite{Billingsley2}, page 363) such that
for all $n\geq m$ and all $a,b\in\S_{k}$
$$
\vert\mu(\xi_{n}^{(k)}=b\vert\xi_{0}^{(k)}=a)-\mu(\xi_{n}^{(k)}=b)\vert<\varepsilon.
$$
Since $\mu(\xi_{n}^{(k)}=x)>0$ for all $x\in\S_{k}$, we can choose
$m=m_{k}$ so that for a given $\varepsilon_{k}>0$
\begin{equation}\label{majoration-fondamentale}
\bigg\vert\frac{\mu(\xi_{0}^{(k)}=a,\xi_{n}^{(k)}=b)}{\mu(\xi_{0}^{(k)}=a)\mu(\xi_{n}^{(k)}=b)}-1\bigg\vert<\varepsilon_{k}
\end{equation}
for all $n\geq m_{k}$ and $a,b\in\S_{k}$.\\
Let $A\in\sigma(\xi_{i}^{(k)},\,i\leq 0\}$ and
$B\in\sigma(\xi_{i}^{(k)},\,i\geq m_{k})$. Suppose first that $A$
and $B$ are elementary cylinders, i.e.
$A=\{\xi_{-s}^{(k)}=e_{-s},...,\xi_{0}^{(k)}=e_{0}\}$ and
$B=\{\xi_{m_{k}}^{(k)}=e_{m_{k}},...,\xi_{t}^{(k)}=e_{t}\}$ with
$e_{0},e_{m_{k}}\in\S_{k}$. Then we derive
$$
\vert\mu(A\cap B)-\mu(A)\mu(B)\vert
=\mu(A)\mu(B)\bigg\vert\frac{\mu(\xi_{m_{k}}^{(k)}=e_{m_{k}},\xi_{0}^{(k)}=e_{0})}{\mu(\xi_{m_{k}}^{(k)}=e_{m_{k}})\mu(\xi_{0}^{(k)}=e_{0})}-1\bigg\vert.
$$
From ($\ref{majoration-fondamentale}$) it follows
\begin{equation}\label{inegalite-cylindres}
\vert\mu(A\cap
B)-\mu(A)\mu(B)\vert\leq\varepsilon_{k}\mu(A)\mu(B).
\end{equation}
Let us suppose that
$p_{1},...,p_{k-1},p_{k}^{'},\,N_{1},...,N_{k}$ have been chosen,
we choose the $\varepsilon_{k}>0$ and an appropriate $m_{k}$ (such
that $\varepsilon_{k}\downarrow 0$ and $m_{k}\uparrow\infty$).
Notice that for the Markov chain $(\xi_{i})$ the states
$T_{k}^{j}F_{l}=T^{j}F_{l},\,j=0,...,N_{l}-1,\,l=1,...,k-1$ remain
the same when we turn to the $(k+1)$-st step. If $\delta_{k}>0$ is
small enough the probabilities of transitions from
$T^{N_{l}-1}F_{l},\,l=1,...,k$, will be almost the same (as close
as we wish) as the probabilities of transitions from
$T_{k}^{N_{l}-1}F_{l},\,l=1,...,k-1,\,T_{k}^{N_{k}-1}F_{k}^{'}$ to
$F_{l},\,l=1,...,k-1,\,F_{k}^{'}$. Therefore, if we set
$\delta_{k}>0$ small enough, for each elementary cylinder $A$ in 
$\sigma(\xi_{i},\,i\leq 0)$ and for each elementary
cylinder $B$ in $\sigma(\xi_{i},\,i\geq m_{k})$ such that the 0-st
coordinate of $A$ and the $m_{k}$-coordinate $B$ are in
$\Delta_{k}=\cup_{l=1}^{k}\cup_{j=0}^{N_{l}-1}T^{j}F_{l}$, we can 
deduce (using ($\ref{inegalite-cylindres}$))
\begin{equation}\label{inegalite-cylindres2}
\vert\mu(A\cap B)-\mu(A)\mu(B)\vert\leq
\varepsilon_{k}\mu(A)\mu(B).
\end{equation}
Now if $A$ and $B$ are finite disjoint unions of such elementary
cylinders, i.e. $A=\cup_{i\in I}A_{i}$ and $B=\cup_{j\in J}B_{j}$,
then
\begin{align*}
\vert\mu(A\cap B)-\mu(A)\mu(B)\vert&\leq\sum_{(i,j)\in I\times
J}\vert\mu(A_{i}\cap B_{j})-\mu(A_{i})\mu(B_{j})\vert\\
&\leq\varepsilon_{k}\sum_{(i,j)\in I\times
J}\mu(A_{i})\mu(B_{j})\\
&=\varepsilon_{k}\mu(A)\mu(B).
\end{align*}
Consequently the inequality ($\ref{inegalite-cylindres2}$) still
holds for any $A\in\sigma(\xi_{i},\,i\leq 0)$ with 0-st coordinate
in $\Delta_{k}$ and
$B\in\sigma(\xi_{i},\,i\geq m_{k})$ with $m_{k}$-st coordinate in $\Delta_{k}$ (in fact, such measurable sets can be
approximated by finite disjoint unions of elementary cylinders).\\
If we choose $\delta_{k}>0$ small enough, the measure of
$\Delta_{k}$ can be made arbitrarily close to 1, hence we get for any $A$ in $\sigma(\xi_{i},\,i\leq 0)$ and any $B$ in $\sigma(\xi_{i},\,i\geq m_{k})$
\begin{equation}\label{inegalite-importante}
\vert\mu(A\cap
B)-\mu(A)\mu(B)\vert<\varepsilon_{k}\mu(A)\mu(B)+6\mu(\Delta_{k}^{c})\leq
7\varepsilon_{k}\converge{10}{k}{+\infty}{}0.
\end{equation}
Since in the definition
($\ref{def-f}$) of the function $f$, the i.i.d. sequence $(g\circ
T^{k})_{k\geq 0}$ is independent of the set $A$ then the process
$(f\circ T^{k})_{k\geq 0}$ is strongly mixing if and only if the
process $(\ind{A^{c}}\circ T^{k})_{k\geq 0}$ is. Actually, we can
notice that there exists a measurable function $h$ such that for
any integer $i\geq 0$
$$
\ind{A^{c}}\circ T^{i}=h(\xi_{i})
$$
where $\xi$ is the Markov chain with state space $\S$ defined
above. Using ($\ref{inegalite-importante}$) and noting that the
$\sigma$-algebras $\sigma(\xi_{i},\,i\leq 0)$ and
$\sigma(\xi_{i},\,i\leq m_{k})$ contain respectively the
$\sigma$-algebras $\sigma(h(\xi_{i}),\,i\leq 0)$ and
$\sigma(h(\xi_{i}),\,i\leq m_{k})$ we obtain the strong mixing
property.$\qquad\qquad\Box$
\begin{Rmq} {\em The process $(f\circ T^i)$, which gives the counterexample of Theorem 3 can in fact
be constructed in any dynamical system of positive entropy.
Like in the proof of Theorem 1, it is sufficient to show that it can be constructed
in any Bernoulli shift.\\
In the proof of Theorem 3 we constructed a partition $\S = \{T^jF_k,
0\leq j\leq N_k-1, k=1,2,\dots\}$ which generates the $\sigma$-field $\A$. If we let the sequence of $\varepsilon_k$ converge to 0 sufficiently fast we can (using ($\ref{inegalite-cylindres}$)) guarantee that for any
$\varepsilon>0$ there exists $N$ such that for all $m$ the partition $\bigvee_{-m}^{0}T^{i}\S$ is $\varepsilon$-independent of $\bigvee_{N}^{N+m}T^{i}\S$
(such dynamical system is called ``weakly Bernoulli").
This implies that the dynamical system is Bernoulli (cf. \cite{Rudolph-Weiss}, \cite{Shields}).
Unfortunately, we did not find any article or book where the implication is
presented for the case of countably infinite partitions (there are numerous references for the case of finite partitions, e.g. \cite{Shields}).}
\end{Rmq}
\textbf{Acknowledgements}\\
\\
We thank Professors M. Denker, D. Rudolph and B. Weiss and the anonymous referee for valuable suggestions which improved the presentation of this paper.
% ------------------------------------------------------------------------
\bibliographystyle{Paul1}
\bibliography{xbib}
\vspace{1cm}
Mohamed EL MACHKOURI, Dalibor VOLN\'Y\\
Laboratoire de Math\'ematiques Rapha\"el Salem\\
UMR 6085, Universit\'e de Rouen\\
Site Colbert\\
76821 Mont-Saint-Aignan, France\\
mohamed.elmachkouri@univ-rouen.fr\\
dalibor.volny@univ-rouen.fr
\end{document}